\documentclass[12pt]{amsart}
\usepackage{amsfonts, amsmath}
\usepackage{amssymb}
\usepackage{graphicx}
\usepackage{epstopdf}
\usepackage{epsfig}
\usepackage{psfrag}
\usepackage{color}

\theoremstyle{plain}
\newtheorem{theorem}{Theorem}[section]
\newtheorem{corollary}[theorem]{Corollary}
\newtheorem{lemma}[theorem]{Lemma}

\theoremstyle{definition}
\newtheorem{definition}[theorem]{Definition}

\newtheorem{notation}[theorem]{Notation}

\newtheorem{problem}[theorem]{Problem}

\newtheorem{remark}[theorem]{Remark}

\begin{document}
\title[Stabilization, amalgamation, and curves of intersection]{Stabilization, amalgamation, and curves of intersection of Heegaard splittings} 

\numberwithin{figure}{section}
\newcommand{\hs}{$V \cup_S W$ }
\newcommand{\pq}{$P \cup_{\Sigma} Q$ }
\newcommand{\hhs}{$V' \cup_{S'} W'$ }
\newcommand{\ppq}{$P' \cup_{\Sigma'} Q'$ }
\psfrag{M}{$M$}
\psfrag{F}{$F$}
\psfrag{V}{$V$}
\psfrag{W}{$W$}
\psfrag{V'}{$V'$}
\psfrag{W'}{$W'$}
\psfrag{S}{$S$}
\psfrag{X}{$H$}
\psfrag{F0}{$F_0$}
\psfrag{g}{$\gamma$}
\definecolor{medgray}{gray}{0.30}
\psfrag{D}{\color{medgray}{$D$}}
\psfrag{F4}{\color{red}{$F$}}
\psfrag{S4}{\color{cyan}{$S$}}
\psfrag{x1}{$X$}
\psfrag{y1}{$Y$}
\psfrag{x2}{\tiny $X$}
\psfrag{y2}{\tiny $Y$}

\author{Ryan Derby-Talbot}
\address{Department of Mathematics, The American University in Cairo}
\email{rdt@aucegypt.edu}

\begin{abstract}
We address a special case of the Stabilization Problem for Heegaard splittings, establishing an upper bound on the number of stabilizations required to make a Heegaard splitting of a Haken 3-manifold isotopic to an amalgamation along an essential surface. As a consequence we show that for any positive integer $n$ there are 3-manifolds containing an essential torus and a Heegaard splitting such that the torus and splitting surface must intersect in at least $n$ simple closed curves. These give the first examples of lower bounds on the minimum number of curves of intersection between an essential surface and a Heegaard surface that are greater than one. 
\end{abstract}


\maketitle

\section{Introduction}

Two types of surfaces, Heegaard surfaces and incompressible surfaces (especially essential surfaces), have been seen to be particularly useful in studying 3-manifolds. Interestingly, these two kinds of surfaces have opposite compressibility properties: an incompressible surface admits no compressing disks, while a Heegaard surface admits an entire system of compressing disks on both sides of the surface. Schultens observed a nice relationship between these two types of surfaces, defining what is known as an amalgamation \cite{Schultens}. In brief, an amalgamation is a Heegaard splitting created by a kind of generalized connected sum: if a surface $F$ cuts a 3-manifold $M$ into submanifolds $X$ and $Y$, then an amalgamation is a Heegaard splitting of $M$ obtained from Heegaard splittings of $X$ and $Y$  (see Section~\ref{sec:amalgamation} for the precise definition). We will assume that the surface $F$ has non-zero genus, since the case that $F$ is a sphere is somewhat unique and is discussed in \cite{Haken}.

In this paper, we investigate a special case of the Stabilization Problem for Heegaard splittings (see Problems~\ref{prob:StabilizationProblem} and \ref{prob:Stabilization_Amalgamation_Problem}), namely: Determine the number of stabilizations required for a Heegaard splitting of a 3-manifold to be isotopic to an amalgamation along an essential surface. Upon considering this problem, we have the following result: 

\begin{theorem}
\label{thm:main}
Let \hs be a Heegaard splitting of a 3-manifold $M$ such that $S$ intersects a mutually separating 
essential surface $F$ in $k$ simple closed curves. Then $V \cup_S W$ is isotopic to an amalgamation along $F$ after at most $$k - \chi(F)$$ stabilizations.
\end{theorem}

Note that the curves of intersection between $F$ and $S$ need not be essential. This bound may be slightly improved (but not stated as cleanly) by considering the number of components of $F$ cut along $S$ instead of the number of curves of intersection between the surfaces (see Corollary~\ref{cor:bound}). As the proof of Theorem~\ref{thm:main} relies on local techniques near the surface $F$ and does not appeal to the global topology of $M$, there are likely many examples of 3-manifolds where the bound given is not best possible. Theorem~\ref{thm:main} relies heavily on Theorem~\ref{theamalgamationlemma}, which establishes conditions for when a Heegaard splitting is recognizable as an amalgamation along an essential surface. 

The assumption that $F$ is {\em mutually separating} (see Definition~\ref{def:mutuallyseparating}) guarantees that $F$ has the appropriate separating properties to be able to form an amalgamation along $F$. A surface that is not mutually separating can be made so by adding parallel copies of some of its components, in which case we can apply Theorem~\ref{thm:main} to this modified surface. In Section~\ref{sec:nonseparating} we discuss this issue in greater detail. 

An interesting application of Theorem~\ref{thm:main} is that it can be applied to 3-manifolds that have ``degeneration of Heegaard genus'' to give a lower bound on the number of curves of intersection between a minimal genus Heegaard splitting surface and an essential surface (see Theorem~\ref{thm:degeneration_of_Heegaard_genus}). Schultens and Weidmann have given examples in  \cite{SchultensWeidmann} where degeneration of Heegaard genus can be arbitrarily large. Combining these results leads to the following theorem:

\begin{theorem}
\label{thm:lowerbound}
For every positive integer $n$, there exist 3-manifolds $M_n$ containing an essential torus $T_n$ and a Heegaard splitting $V_n \cup_{S_n} W_n$ such that the minimum number of simple closed curves of intersection between $T_n$ and $S_n$ is at least $n$.
\end{theorem}

While there have been previous results establishing {\em upper} bounds on the minimum number of curves of intersection between a Heegaard surface and an essential surface (see Section~\ref{sec:bounds}), these are the first examples giving non-trivial {\em lower} bounds on the minimum number of curves of intersection. 

This paper is organized as follows. In Section~\ref{sec:Heegaard_splittings} we provide the background and basic definitions concerning Heegaard splittings, stabilization and amalgamation. In Section~\ref{sec:stabilization} we prove Theorem~\ref{theamalgamationlemma} which allows us to determine when a Heegaard splitting is an amalgamation, and use it to prove Theorem~\ref{thm:main}. In Section~\ref{sec:bounds} we discuss bounds on the number curves of intersection of Heegaard splittings and essential surfaces, establish lower bounds when a 3-manifold has degeneration of Heegaard genus (Theorem~\ref{thm:degeneration_of_Heegaard_genus}), and finally prove Theorem~\ref{thm:lowerbound}. 

\section{Heegaard splittings}
\label{sec:Heegaard_splittings}

\subsection{Basic definitions}

In this paper, $M$ denotes a compact, orientable 3-manifold. All surfaces in $M$ are assumed to be orientable and embedded.

\begin{definition}
Let $S$ be a closed, orientable surface. A {\em compression body} $V$ is a 3-manifold obtained from $S \times I$ by attaching $2$-handles to $S \times \{0\}$ and capping off any resulting $2$-sphere components with 3-balls. We denote $ \partial_+ V = S \times \{1\} = S$, and $\partial_- V = \partial V - \partial_+V$. In the case that $\partial_- V = \emptyset$, $V$ is called a {\em handlebody}.
\end{definition}

Dually, $V$ can be obtained by taking a closed orientable surface $\widetilde{S}$ and attaching 1-handles to $\widetilde{S} \times I$ along $\widetilde{S} \times \{1\}$ (or to a 3-ball along the boundary if $V$ is a handlebody). In this case $\partial_- V = \widetilde{S} \times \{0\}$ and $\partial_+ V = \partial V - \partial_- V$. In either case, the {\em genus} of $V$ is the genus of the surface $\partial_+ V$. 

\definecolor{light-gray}{gray}{0.5}
\psfrag{-}{\color{light-gray}{$\partial_- V$}}
\psfrag{+}{$\partial_+ V$}

\begin{figure}[h] 
   \centering
   \includegraphics[width=3.5in]{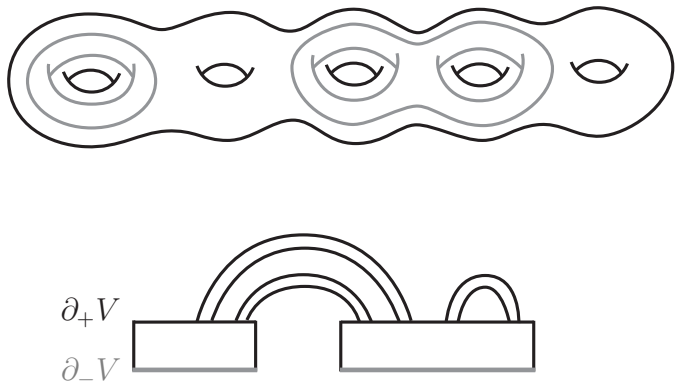} 
   \caption{A genus 5 compression body (above) and a schematic for it (below) illustrating the attachment of 1-handles to $\partial_- V \times I$.}
   \label{fig:compressionbody}
\end{figure}

\begin{definition}
A {\em spine} of a compression body $V$ is a graph $\sigma$ embedded in $V$ with $\partial \sigma \cap V = \sigma \cap \partial_- V$ such that $V$ deformation retracts onto $\sigma \cup \partial_- V$. 
\end{definition}

\begin{figure}[h] 
   \centering
   \includegraphics[width=3in]{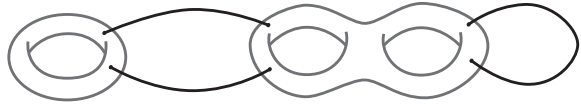} 
   \caption{A spine of the compression body in Figure~\ref{fig:compressionbody}.}
   \label{fig:spine}
\end{figure}

Any two spines of a compression body are equivalent up to isotopy and edge slides.

\begin{definition}
A {\em Heegaard splitting} $V \cup_S W$ of $M$ is a decomposition of $M$ into two compression bodies $V$ and $W$ of the same genus such that $M$ is obtained from $V$ and $W$ by identifying $\partial_+ V$ with $\partial_+ W$ via some homeomorphism. We denote the surface $S = \partial_+ V = \partial_+ W$ in $M$ as the {\em Heegaard surface} or {\em splitting surface} of $V \cup_S W$. The {\em genus} of $V \cup_S W$ is the genus of $S$.
\end{definition}

\psfrag{r}{$\partial M$}
\psfrag{s}{\color{cyan}{$S$}}

\begin{figure}[h] 
   \centering
   \includegraphics[width=2.5in]{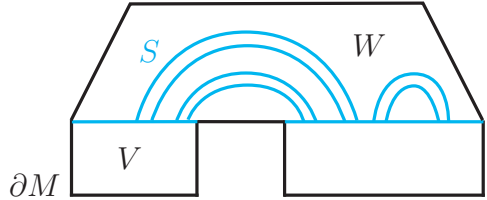} 
   \caption{A schematic of a Heegaard splitting \hs of $M$. In this figure, $V$ is a compression body (the same one as in Figure~\ref{fig:compressionbody}) with two components of $\partial_- V$ equalling $\partial M$, and $W$ is a handlebody. The Heegaard surface $S$ is given in blue.}
   \label{fig:heegaardsplitting}
\end{figure}

Two Heegaard splittings $V \cup_S W$ and $P \cup_{\Sigma} Q$ of $M$ are {\em isotopic} if there is an isotopy of $M$ taking $V$ to $P$. 

\subsection{Stabilization}

A classic result of Moise \cite{Moise} implies that every 3-manifold $M$ has a Heegaard splitting. It follows that $M$ admits infinitely many Heegaard splittings (up to isotopy) via the following construction. 

\begin{definition}
Let $V \cup_S W$ be a Heegaard splitting of $M$. Add a 1-handle $H$ to $V$ along $\partial_+ V$ such that its core is isotopic in $M$ to an arc on $\partial_+ V$ ({\em i.e.}~the core is unknotted). Then $V' = V \cup H$ and $W' = \overline{W - H}$ are compression bodies. The resulting Heegaard splitting $V' \cup_{S'} W'$ of $M$ is obtained by a {\em stabilization} of $V \cup_S W$. 
\end{definition}

\begin{figure}[h]
   \centering
   \includegraphics[width=3in]{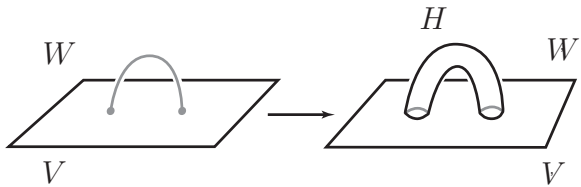} 
   \caption{A stabilization of $V \cup_S W$.}
   \label{fig:Stabilization}
\end{figure}

Note that the genus of $V' \cup_{S'} W'$ is one greater than the genus of $V \cup_S W$. We sometimes refer to the splitting $V' \cup_{S'} W'$ itself as a {\em stabilization} of $V \cup_S W$, or as a {\em stabilized} Heegaard splitting. It is a fact that $V' \cup_{S'} W'$ is stabilized if and only if there exist essential disks $D \subset V'$ and $E \subset W'$ that intersect in a single point (see {\em e.g.}~\cite{Scharlemann}). Note that a stabilization of a Heegaard splitting is unique in $M$ up to isotopy.

It is a classic theorem of Reidemeister \cite{Reidemeister} and Singer \cite{Singer} that any two Heegaard splittings of $M$ can be made isotopic after a sufficient number of stabilizations of each splitting (assuming the Heegaard splittings partition $\partial M$ in the same way). The question remains, however, as to the number of stabilizations needed to achieve isotopy. This is called the {Stabilization Problem}. 

\begin{problem}[The Stabilization Problem]
\label{prob:StabilizationProblem}
Given two Heegaard splittings of a 3-manifold $M$, determine the minimum number of stabilizations required to make the two splittings isotopic. 
\end{problem}

Several examples are known where only one stabilization (of the larger genus splitting) is needed to achieve isotopy (see {\em e.g.}~\cite{DerbyTalbot2006}, \cite{Hagiwara}, \cite{SchultensSFS}, and \cite{Sedgwick}). Recently, Bachman \cite{Bachman} and independently Hass, Thompson and Thurston \cite{Hass} as well as Johnson \cite{Johnson2}, \cite{Johnson} have shown that the necessary number of stabilizations can be much greater than one, in fact as large as $g$, where $g$ denotes the genera of the initial Heegaard splittings.

In general, establishing upper bounds on the number of stabilizations required to make two Heegaard splittings isotopic is a difficult problem. Rubinstein and Scharlemann have shown that if $M$ is non-Haken, then two Heegaard splittings of $M$ of genus $g$ and $g'$, respectively, with $g \geq g'$ are isotopic after at most $7g + 5g' - 9$ stabilizations of the larger genus splitting \cite{RubinsteinScharlemann}. In previous work \cite{DerbyTalbot2007}, the author showed that under mild assumptions, two Heegaard splittings of genus $g$ obtained by Dehn twisting along a JSJ torus in $M$ are isotopic after at most $4g-4$ stabilizations.  In neither case have these bounds been shown to be sharp. Establishing better bounds is a rich area for future research.

\subsection{Amalgamation}
\label{sec:amalgamation}

\begin{definition}
A {\em compressing disk} for a properly embedded surface $F$ in $M$ is an embedded disk $D$ such that $D \cap F =  \partial D$, and $\partial D$ does not bound a disk in $F$. A surface $F$ is called {\em incompressible} if $F$ admits no compressing disks. A surface $F$ is called {\em essential} if $F$ is incompressible and no component of $F$ is boundary parallel. 
\end{definition}

We now present a technique originally due to Schultens \cite{Schultens} of constructing Heegaard splittings of $M$ from Heegaard splittings of components of $M$ obtained by cutting along some surface $F$. Almost always this surface is taken to be incompressible. Moreover, if (a component of) $F$ is boundary parallel, then {\em any} Heegaard splitting \hs of $M$ is a ``trivial'' amalgamation of itself and a ``Type I'' splitting of $F \times I$ \cite{ScharlemannThompson} along $F$, usually not considered to be an amalgamation. In light of this observation we will henceforth take $F$ to be essential. 

Although the procedure for defining an amalgamation is straightforward, one must take care in describing how the component Heegaard splittings are attached together so that the resulting splitting surface remains separating in the manifold. See {\em e.g.}~\cite{BachmanDT} or \cite{Lackenby} for descriptions of the construction in simpler cases.

\begin{definition}
\label{def:mutuallyseparating}
A closed, orientable surface $F$ in a 3-manifold $M$ is called {\em mutually separating} if $M$ cut along $F$ consists of two (possibly disconnected) 3-manifolds $X$ and $Y$ such that every neighborhood of each component of $F$ intersects both $X$ and $Y$. 
\end{definition}

In essence, a surface $F$ is mutually separating in $M$ if we can color the components of $M$ cut along $F$ alternately in black and white.

\begin{figure}[h] 
\centering
\includegraphics[width=5.5in]{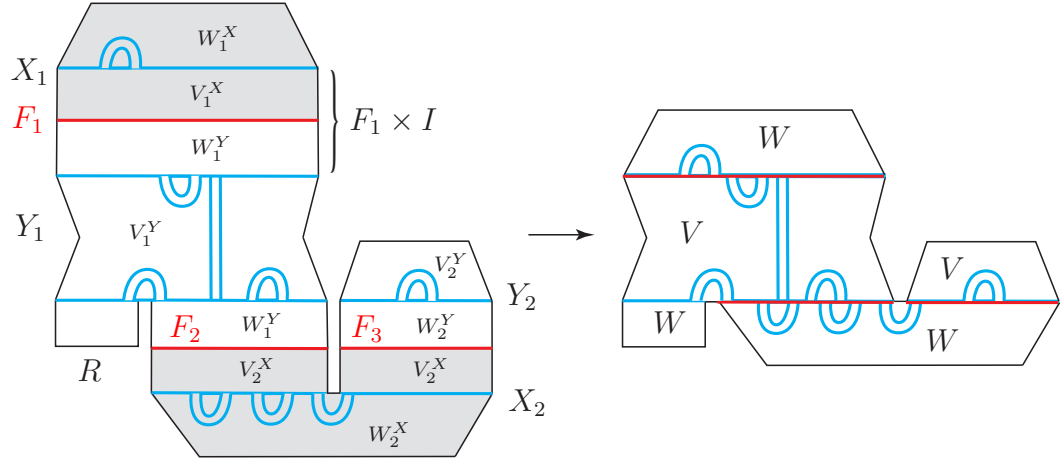}
   \caption{Forming an amalgamation along $F$: in this schematic the surface $F$ consists of three components $F_1$, $F_2$, $F_3$ (in red), cutting the manifold into submanifolds $X_1$, $X_2$ (shaded) and $Y_1$, $Y_2$. The Heegaard surfaces $S_i^X$ and $S_i^Y$, $i = 1,2$, are colored in blue. The surface $R$ denotes a component of $\partial M$ contained in $\partial_- W_1^Y$.}
   \label{fig:amalgamation}
\end{figure}

\begin{definition}
\label{def:amalgamation}
Let $F$ be a mutually separating surface in $M$ such that $M$ cut along $F$ equals $X$ and $Y$ as above. Let $X_1, \ldots, X_m$ and $Y_1, \ldots, Y_n$ be the components of $X$ and $Y$, respectively. For each $i, 1 \leq i \leq m$, let $V_i^X \cup_{S_i^X} W_i^X$ be a Heegaard splitting of $X_i$ such that $F \cap X_i \subset \partial_- V_i^X$. Similarly for each $j$, $1 \leq j \leq n$, let $V_j^Y \cup_{S_j^Y} W_j^Y$ be a Heegaard splitting of $Y_j$ such that $F \cap Y_j \subset \partial_- W_j^Y$. The surface $F$ has a product neighborhood $N(F)$ such that in each $X_i$, $N(F) \cap X_i \subset V_i^X$ and $V_i^X - N(F)$ consists of 1-handles and possibly a component homeomorphic to $Q \times I$ where $Q$ is a subset of the components of $\partial X_i - F$. Similarly, in each $Y_j$, $N(F) \cap Y_j \subset W_j^Y$ and $W_j^Y - N(F)$ consists of 1-handles and possibly a component homeomorphic to $R \times I$ where $R$ is a subset of the components of $\partial Y_j - F$. In each component $\widetilde{F} \times I$ of $N(F)$, identify $\widetilde{F} \times \{t\}$ with $\widetilde{F} \times \{ \frac{1}{2} \}$ for $t \in I$, so that the ends of any 1-handles in $V_i^X$ or $W_j^Y$ meeting $\widetilde F \times I$ are taken to be disjoint on $\widetilde{F} \times \{\frac{1}{2}\}$. The resulting manifold is homeomorphic to $M$, and 
\[ V = \bigcup_{i=1}^m \left(V_i^X - N(F)\right) \ \cup \ \bigcup_{j=1}^n V_j^Y \]
\[ W = \bigcup_{j=1}^n \left(W_j^Y - N(F)\right) \ \cup \ \bigcup_{i=1}^m W_i^X \] 
are compression bodies. The resulting Heegaard splitting $V \cup_S W$ of $M$ is called an {\em amalgamation along $F$}. (See Figure~\ref{fig:amalgamation}.) 
\end{definition}

\psfrag{X1}{$X_1$}
\psfrag{X2}{$X_2$}
\psfrag{Y1}{$Y_1$}
\psfrag{Y2}{$Y_2$}
\psfrag{F1}{\color{red}{$F_1$}}
\psfrag{F2}{\color{red}{$F_2$}}
\psfrag{F3}{\color{red}{$F_3$}}
\psfrag{V1}{\tiny$V_1^X$}
\psfrag{V2}{\tiny$V_2^X$}
\psfrag{W1}{\tiny$W_1^X$}
\psfrag{W2}{\tiny$W_2^X$}
\psfrag{U1}{\tiny$V_1^Y$}
\psfrag{U2}{\tiny$V_2^Y$}
\psfrag{R1}{\tiny$W_1^Y$}
\psfrag{R2}{\tiny$W_2^Y$}
\psfrag{p}{$R$}
\psfrag{times}{$F_1 \times I$}

The simplest case of the above definition is when $F$ is connected and separating. Then an amalgamation along $F$ is obtained from only two Heegaard splittings, one from each of the components of $M$ cut along $F$. If $F$ is non-separating and connected, then $F$ is clearly not mutually separating and so the definition of amalgamation does not make sense. We will consider how to deal with this situation in Section~\ref{sec:nonseparating}. For convenience of discussion, assume for the rest of this section that $F$ is mutually separating. 

\begin{remark}
\label{rmk:F_minus_open_disks}
Notice that if \hs is an amalgamation along $F$, then $S$ can be isotoped to intersect $F$ in $F - \{open \ disks\}$, as indicated by the right hand side of Figure~\ref{fig:amalgamation} (the open disks being where the 1-handles attach). The converse is also true: if $S\cap F = F - \{open \ disks\}$, then \hs can be untelescoped into Heegaard splittings of the components of $M$ cut along $F$, implying that \hs is an amalgamation along $F$ by definition.
\end{remark}

\subsection{The Stabilization-Amalgamation Problem}
\label{sec:stabilizationamalgamationproblem}

Recent work has shown that in some sense (low genus) Heegaard splittings of ``generic'' 3-manifolds are amalgamations (see {\em e.g.}~\cite{Bachman2008}, \cite{BSS}, \cite{Lackenby}, \cite{Li} and \cite{Souto}). Not every Heegaard splitting, however, is an amalgamation of Heegaard splittings along some essential surface. For example, so-called strongly irreducible Heegaard splittings, first defined in \cite{CassonGordon} by Casson and Gordon, are not amalgamations. There are many examples of 3-manifolds admitting strongly irreducible Heegaard splittings. Moreover, a Heegaard splitting may be an amalgamation along one surface but not necessarily along another. In light of these observations, one can ask the following special case of the Stabilization Problem:

\begin{problem}[The Stabilization-Amalgamation Problem]
\label{prob:Stabilization_Amalgamation_Problem}
Determine the minimum number of stabilizations required to make a Heegaard splitting of $M$ and an amalgamation along $F$ isotopic. 
\end{problem}

A solution to the Stabilization-Amalgamation Problem would provide a strategy for addressing the Stabilization Problem for Haken 3-manifolds in the following way. Let \hs and \pq be Heegaard splittings of a 3-manifold $M$ containing an essential surface $F$. Suppose that \hs and \pq require $s$ stabilizations to become isotopic to (stabilized) amalgamations along $F$. Then, the stabilized splittings \hhs and \ppq are constructed from Heegaard splittings of the components of $M$ cut along $F$ via amalgamation. Moreover, assuming that $V' \cup_{S'} W'$ and $P' \cup_{\Sigma'} Q'$ partition $\partial M$ in the same way, then the splittings of the components of $M$ cut along $F$ forming the amalgamation can also be assumed to partition the boundary components of the components of $M$ cut along $F$ in the same way. If it is known that the maximum number of stabilizations required for two Heegaard splittings of the components of $M$ cut along $F$ with the given boundary partitions to become isotopic is $s'$, then at most an additional $s' - 1$ stabilizations are needed (since the splittings are already stabilized at least once) for the Heegaard splittings in a given component of $M$ cut along $F$ to become isotopic. After being used for this isotopy, these stabilizations can then be passed to the next component where the process is repeated, and so on. Thus \hs and \pq are isotopic after at most $s + s' -1$ stabilizations. 

This strategy is utilized in \cite{DerbyTalbot2006}, where it is shown that $s = s' = 1$ for strongly irreducible Heegaard splittings of totally oriented graph manifolds (assuming that one of the splittings has genus at least as large as a minimal genus amalgamation along the JSJ tori of $M$). More specifically, if $\Theta$ denotes the canonical system of JSJ tori of $M$, then it is shown that a strongly irreducible Heegaard splitting \hs of $M$ and an amalgamation along $\Theta$, $P \cup_{\Sigma} Q$, become isotopic after at most one stabilization of the larger genus splitting. This is done by first isotoping the stabilized splitting $V' \cup_{S'} W'$ to be an amalgamation along $\Theta$ (this relies heavily on exploiting the underlying structure of the graph manifold). Then $V' \cup_{S'} W'$ and $P' \cup_{\Sigma'} Q'$ can be untelescoped into Heegaard splittings of the Seifert fibered components of $M$. The splittings in each component can be made isotopic by appealing to the main result of \cite{SchultensSFS} which says that any two stabilized Heegaard splittings of a Seifert fibered space are isotopic (assuming that they have the same genus, which can be assumed to be true in this case). The stabilization(s) of each Heegaard splitting can be passed from Seifert fibered component to Seifert fibered component to obtain this isotopy. Thus this procedure can be applied to two strongly irreducible Heegaard splittings to obtain the aforementioned result. 

In light of this strategy, it becomes of interest to study the Stabilization-Amalgamation Problem for Haken 3-manifolds. This is our undertaking in the next sections. 

\section{Stabilization versus Amalgamation}
\label{sec:stabilization}

In this section we establish a relationship between Heegaard splittings that are stabilizations and those that are amalgamations. 

\subsection{Compressing structures on surfaces}

\begin{notation}
For a surface $F$, let $g(F)$ denote the genus of $F$. If $G$ is an $n$-manifold, then let $|G|$ denote the number of components of $G$. We will be primarily concerned with the case that $G$ is a 1- or 2-manifold. 
\end{notation}

\begin{definition}
Let $F_0$ be a surface with boundary. An embedded arc $\gamma$ in $F_0$ is a {\em compressing arc} (also called an {\em essential arc}) if $\gamma \cap \partial F_0 = \partial \gamma$ and $\gamma$ does not cut off a disk from $F_0$.
\end{definition}

Note that a disk admits no compressing arcs. 

\begin{definition}
Let $F_0$ be a surface with boundary. A {\em complete system of compressing arcs} $\Gamma$ for $F_0$ is a disjoint union of compressing arcs in $F_0$ such that $F_0$ cut along $\Gamma$ is a disjoint union of $|F_0|$ disks. 
\end{definition}

\begin{remark}
\label{numberofcompressingarcs}
Note that if $F_0$ is disconnected such that every component has boundary, then a complete system of compressing arcs for each component $\widetilde F_0$ of $F_0$ has $$2g(\widetilde F_0) + |\partial \widetilde F_0 | - 1 = 1 - \chi(\widetilde F_0)$$ components. (Note that if a component $\widetilde{F}_0$ is a disk, then it admits no compressing arcs and so the above number is zero.) Since Euler characteristic is additive, a complete system of compressing arcs for $F_0$ has $$|F_0| - \chi(F_0)$$ components. 
\end{remark}

\begin{figure}[h]
   \centering
   \includegraphics[width=1.2in]{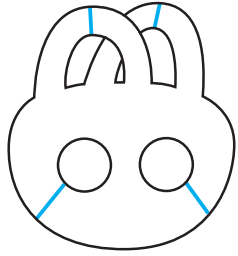} 
   \caption{A complete system of compressing arcs for a three-times punctured torus.}
   \label{fig:compressingarcs}
\end{figure}

\begin{figure}[h] 
   \centering
   \includegraphics[width=2in]{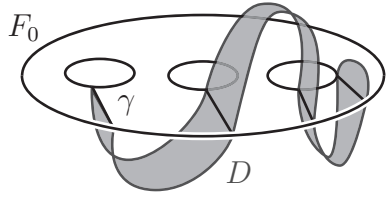} 
   \caption{A weak boundary compressing disk $D$ for $F_0$ based at $\gamma$.}
   \label{fig:weakboundarycompressingdisk}
\end{figure}

\begin{definition}
Let $F_0$ be a surface with boundary properly embedded in a 3-manifold $M$. An embedded disk $D$ is a {\em boundary compressing disk} (or $\partial$-{\em compressing disk}) for $F_0$ if $\partial D = \gamma \cup \delta$ where $D \cap F_0 = \partial D \cap F_0 = \gamma$ is a compressing arc of $F_0$, and $D \cap \partial M = \partial D \cap \partial M = \delta$. In this case we say that $D$ is {\em based at} $\gamma$. 
\medskip

\noindent An embedded disk $D$ is a {\em weak boundary compressing disk} (or {\em weak $\partial$-compressing disk}) for $F_0$ if $\partial D= \gamma \cup \delta$ where $\gamma$ is a compressing arc in $F_0$, $D \cap \partial M = \partial D \cap \partial M = \delta$, and $D-\gamma$ intersects $F_0$ transversely in compressing arcs of $F_0$ that are not parallel in $F_0$ to $\gamma$ (see Figure~\ref{fig:weakboundarycompressingdisk}). As before, we say $D$ is {\em based at} $\gamma$.
\end{definition}

Note that a $\partial$-compressing disk based at $\gamma$ is a weak $\partial$-compressing disk with no additional curves of intersection with $F_0$ (other than $\gamma$). 

\begin{definition}
Let $F_0$ be a surface properly embedded in a 3-manifold $M$, and assume that every component of $F_0$ has boundary. Then a union $\Delta$ of disjoint (weak) $\partial$-compressing disks for $F_0$ is a {\em complete system of (weak) $\partial$-compressing disks} for $F_0$ if there is a complete system $\Gamma$ of compressing arcs for $F_0$ such that for each component $\gamma$ of $\Gamma$ there is exactly one component of $\Delta$ based at $\gamma$. 

\end{definition}

We now state a lemma that will be useful in the next section. 

\begin{lemma}
\label{weakboundarycompressingdisklemma}
Let $F_0$ be an incompressible surface properly embedded in a compression body $V$ so that every component of $F_0$ has boundary, and assume that $\partial F_0 \subset \partial_+V$. Then $F_0$ admits a complete system of weak $\partial$-compressing disks $\Delta$ in $V$. Moreover, at least one of the components of $\Delta$ is a $\partial$-compressing disk. 
\end{lemma}

\begin{proof}
Let $\mathcal D$ be a complete system of compressing disks for $V$, so that $V$ cut along $\mathcal D$ is a 3-ball if $V$ is a handlebody, or is homeomorphic to $\partial_-V \times I$ if $\partial_- V$ is nonempty. In the latter case, let $\sigma$ be a collection of essential simple closed curves on $\partial_-V$ that cut $\partial_- V$ into $|\partial_-V|$ disks. The curves in $\sigma$ can be chosen so that on each component of $\partial_-V$, the curves are disjoint except for a common basepoint, such as a set of curves that generate a basis for the fundamental group. Now $\sigma \times I$ forms a 2-complex in $V$ that can be assumed to be disjoint from the disk components in $\mathcal D$, and that consists of vertical annuli such that in each component of $\partial_-V \times I$ the corresponding annuli in $\sigma \times I$ all intersect in a single vertical arc. Moreover, $\sigma \times I$ cuts each component of $\partial_- V \times I$ into 3-balls. Thus, if $V$ is a compression body, redefine $\mathcal D$ to be the disjoint union of the complete system of compressing disks along with the 2-complex $\sigma \times I$. 

By a standard innermost disk, outermost arc argument, we can assume that $F_0$ intersects each disk component only in compressing arcs of $F_0$. We claim that the same can be done for each annulus $(curve) \times I$ in $\sigma \times I$. First, eliminate any inessential curves of intersection of $F_0$ with each $(curve) \times I$ by performing a standard innermost disk argument. Now let $\Sigma$ be a component of $\sigma \times I$, thus $\Sigma$ consists of vertical annuli in a component of $\partial_- V \times I$ intersecting in a common vertical arc $a$. Then each component of $F_0 \cap \Sigma$ disjoint from $a$ is is an arc with both endpoints on $\partial_+ V$ in one of the annuli in $\Sigma$, and each component of $F_0 \cap \Sigma$ intersecting $a$ is a 1-complex consisting of arcs and essential loops in the vertical annuli of $\Sigma$. 

In the latter case, assume that a component of $F_0 \cap \Sigma$ consists only of essential loops on the annuli of $\Sigma$ all meeting at a common basepoint on $a$. Let $\widetilde{F}_0$ be the component of $F_0$ meeting $\Sigma$ in this way. Then $\widetilde{F}_0$ cut along $\widetilde{F}_0 \cap \Sigma$ must be a disk as $F_0$ is incompressible. But this implies that $\widetilde{F}_0$ has no boundary components, a contradiction. Thus each component of $F_0 \cap \Sigma$ that intersects $a$ must meet at least one of the annuli in $\Sigma$ in an arc with both endpoints on $\partial_+V$. It is now straightforward to perform an isotopy in a neighborhood of $\Sigma$ that moves an outermost such arc off of $a$. Such an isotopy also turns any essential loops in the same component of $F_0 \cap \Sigma$ as the aforementioned arc into arcs on annuli in $\Sigma$ by removing the intersection point with $a$. Thus by repeating this process we can assume that $F_0$ intersects each of the annuli of $\Sigma$ in arcs that have both boundary components on $\partial_+V$ and are disjoint from $a$. Finally, we can perform a standard outermost arc argument to eliminate any arcs of intersection that are not compressing arcs of $F_0$. Thus $F_0$ intersects $\mathcal D$ only in compressing arcs of $F_0$. 
 
Since each component of $V$ cut along $\mathcal D$ is a 3-ball and since $F_0$ is incompressible, the arcs in $\mathcal D \cap F_0$ must cut $F_0$ into disks. Therefore there is a subset $\Gamma$ of $\mathcal D \cap F_0$ such that $\Gamma$ is a complete system of compressing arcs for $F_0$. Moreover, we can choose this subset to include at least one outermost arc of $\mathcal D \cap F_0$, implying the last conclusion of the lemma. Now, given an arc $\gamma$ in $\Gamma$, suppose $D$ is the component of $\mathcal D$ containing $\gamma$ and let $D_{\gamma}$ be a disk component of $D$ cut along $\gamma$. If $(D_{\gamma} \cap F_0) - \gamma$ contains another arc $\gamma'$ in $F_0$ parallel to $\gamma$, then take an outermost such arc and rename it as the arc $\gamma$ in $\Gamma$. Doing this for all the components of $\Gamma$ yields a complete system of weak $\partial$-compressing disks for $F_0$ in $V$. 
\end{proof}

\subsection{Stabilizing to amalgamation}
\label{sec:stabilizationtoamalgamation}

The following theorem establishes when a Heegaard splitting is an amalgamation along an essential surface in $M$. It is a generalization of Lemma 3.1 in \cite{DerbyTalbot2007} (where $F$ was assumed to be a torus).  

\begin{theorem}
\label{theamalgamationlemma}
Let $M$ be a 3-manifold such that $F$ is an essential mutually separating surface cutting $M$ into $X$ and $Y$. A Heegaard splitting $V \cup_S W$ of $M$ is an amalgamation along $F$ if and only if $S$ is isotopic to a surface intersecting $F$ transversely such that $F \cap V$ admits a complete system of $\partial$-compressing disks in $X$ and $F \cap W$ admits a complete system of $\partial$-compressing disks in $Y$ (or vice versa).
\end{theorem}

\begin{proof}
By Remark~\ref{rmk:F_minus_open_disks}, $V \cup_S W$ is an amalgamation along $F$ if and only if $S$ can be isotoped so that $S \cap F = F - \{open \ disks\}$. If $V \cup_S W$ is an amalgamation along $F$ isotoped to meet $F$ in this manner, then the open disks are where the 1-handles of the corresponding compression bodies on either side of $M$ are attached. Thus we can label each disk with an ``$X$'' or a ``$Y$'', depending on whether or not the 1-handle meeting $F$ at the disk is in $X$ or $Y$. There must be at least one disk each labeled ``$X$'' and ``$Y$'' since $F$ is essential and cannot be completely contained in either $V$ or $W$. 

\begin{figure}[h] 
   \centering
   \includegraphics[width=2.7in]{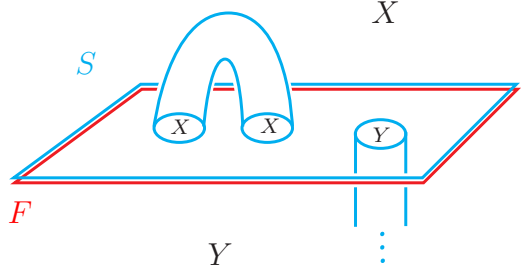} 
   \caption{The Heegaard surface $S$ for an amalgamation along $F$ isotoped to intersect $F$ in $F - \{open \ disks \}$.}
   \label{fig:Amalgamation_and_open_disks}
\end{figure}

Let $\Lambda \subset F$ be a disjoint union of simple closed curves separating $F$ into two (possibly disconnected) subsurfaces $F_X$ and $F_Y$, such that all the disks labeled ``$X$'' are in $F_X$ and all the disks labeled ``$Y$'' are in $F_Y$. Without loss of generality, assume that the ``$X$'' disks are in $V$ and the ``$Y$'' disks are in $W$. Moreover assume that every component $\lambda$ of $\Lambda$ is such that every neighborhood of $\lambda$ intersects both $F_X$ and $F_Y$ ($\lambda$ is ``mutually separating'' in $F$). Then holding $\Lambda$ fixed, isotope $S \cap F_X$ along with any attached 1-handles into $X$ and $S \cap F_Y$ with 1-handles into $Y$ so that after isotopy $F \cap V = F_X$ and $F \cap W = F_Y$. Given complete systems of compressing arcs $\Gamma_X$ of $F_X$ and $\Gamma_Y$ of $F_Y$ disjoint from the original open disks corresponding to the attached 1-handles, the result of pushing $S \cap F_X$ into $X$ and $S \cap F_Y$ into $Y$ gives respective $\partial$-compressing disks in $X$ based at components of $\Gamma_X$ and in $Y$ based at components of $\Gamma_Y$ (see Figure~\ref{fig:amalgamationlemma1}). This proves the forward direction of the lemma.

\psfrag{FX}{\color{red}{$F_X$}}
\psfrag{FY}{\color{red}{$F_Y$}}

\begin{figure}[h] 
   \centering
   \includegraphics[width=2.3in]{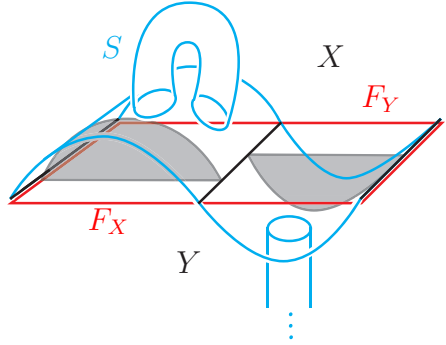} 
   \caption{Isotoping $S$ into $X$ and $Y$ by keeping $\Lambda$ (the black curves) fixed yields a complete system of $\partial$-compressing disks for $F_X$ in $X$ and for $F_Y$ in $Y$.}
   \label{fig:amalgamationlemma1}
\end{figure}

Now assume that $V \cup_S W$ is a Heegaard splitting such that $S$ is isotoped to be transverse to $F$, and such that $F \cap V$ has a complete system $\Delta_{V}$ of $\partial$-compressing disks in $X$, and $F \cap W$ has a complete system $\Delta_{W}$ of $\partial$-compressing disks in $Y$. Let $D$ be a $\partial$-compressing disk of $F \cap V$ (contained in $X$) so that $\partial D = \gamma \cup \delta$, where $\gamma$ is a compressing arc of $F \cap V$ and $\delta \subset S$. Let $N(D)$ be a neighborhood of $D$ such that $N(D) \cap S = \partial N(D) \cap S$ is a product neighborhood of $\delta$ in $S$, and $N(D) \cap F = \partial N(D) \cap F$ is a product neighborhood of $\gamma$ in $F \cap V$. Isotope $S$ ``along $D$'', so that $N(D) \cap S$ is replaced by $\overline{\partial N(D) - (N(D) \cap S)}$ (this operation is known as a {\em boundary compression along $D$}). Doing this isotopy for every component of $\Delta_{V}$ in $X$ and $\Delta_{W}$ in $Y$ (note that the $\partial$-compressing disks are all disjoint by assumption), along with further isotopy of $S$ onto $F$ near $S \cap F$, shows that $S$ can be isotoped to intersect $F$ in $F - \{ open \ disks \}$ (see Figure~\ref{fig:amalgamationlemma2}). Hence $V \cup_S W$ is an amalgamation along $F$.

\psfrag{g}{$\gamma$}
\psfrag{delta}{$\delta$}

\begin{figure}[h] 
   \centering
   \includegraphics[height=2.5in]{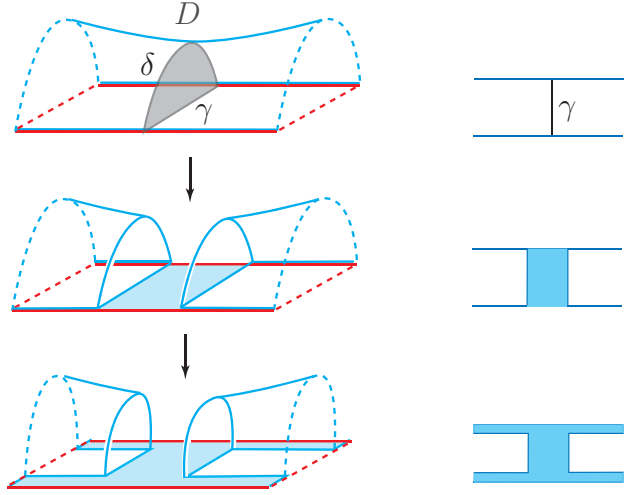} 
   \caption{Isotoping $S$ along $D$ and then near $S \cap F$ so that $S$ intersects $F$ in $F - \{open \ disks\}$. On the right is the corresponding intersection set on $F$ of $S \cap F$ near $\gamma$ for each step.}
   \label{fig:amalgamationlemma2}
\end{figure}

\end{proof}

Theorem~\ref{theamalgamationlemma} yields an answer to the Stabilization-Amalgamation Problem (Problem~\ref{prob:Stabilization_Amalgamation_Problem}) based on a purely local construction near $F$. 

\begin{corollary}
\label{cor:bound}
A Heegaard splitting $V \cup_S W$ is isotopic to an amalgamation along $F$ after at most $$|F \cap V| + |F \cap W| - \chi(F) - 1$$ stabilizations.
\end{corollary}

\begin{proof}
Suppose that $S$ is isotoped to intersect $F$ transversely. We want to apply Theorem~\ref{theamalgamationlemma} to some stabilization $V' \cup_{S'} W'$ of $V \cup_S W$, which means we want to stabilize \hs to obtain complete systems of boundary compressing disks for $F \cap V'$ and $F \cap W'$ in $X$ and $Y$, respectively (or vice versa). By Lemma~\ref{weakboundarycompressingdisklemma}, there exist complete systems of weak $\partial$-compressing disks $\Delta_V$ and $\Delta_W$ for $F \cap V$ and $F \cap W$, respectively. Let $\Gamma_V$ and $\Gamma_W$ be the complete systems of compressing arcs of $F \cap V$ and $F \cap W$, respectively, corresponding to $\Delta_V$ and $\Delta_W$. By the last conclusion of Lemma~\ref{weakboundarycompressingdisklemma}, there must be some component of $\Gamma_V$ or $\Gamma_W$ for which the corresponding weak $\partial$-compressing disk based at that component is actually a $\partial$-compressing disk in $X$ or $Y$. Without loss of generality, assume that this $\partial$-compressing disk is in $X$ and is based at a component of $\Gamma_V$ (note that we cannot also use the lemma to obtain a boundary compressing disk based at a component of $\Gamma_W$ in this way, since it may not necessarily be in $Y$). 

For each remaining component $\gamma$ of $\Gamma_V \cup \Gamma_W$, form $V' \cup_{S'} W'$ by adding a 1-handle to $W$ if $\gamma$ is in $\Gamma_V$ and to $V$ if $\gamma$ is in $\Gamma_W$, so that $\gamma$ is its core. Each added 1-handle is clearly a stabilization, since its cocore and the weak boundary compressing disk component of $\Delta_V$ or $\Delta_W$ based at $\gamma$ intersect in a single point. 

Now, isotope the 1-handles added to $V$ slightly into $Y$ and the 1-handles added to $W$ slightly into $X$. This immediately gives a complete system of $\partial$-compressing disks for $F \cap V'$ in $X$ and also for $F \cap W'$ in $Y$ (see Figure~\ref{fig:stabilization_boundary_compressing_disk}). Thus Theorem~\ref{theamalgamationlemma} applies to show that the stabilized splitting is an amalgamation along $F$. By Remark~\ref{numberofcompressingarcs}, the number of stabilizations used is $|F \cap V| + |F \cap W| - \chi(F) - 1$.

\psfrag{rF}{\color{red}{$F$}}
\psfrag{cS}{\color{cyan}{$S'$}}

\begin{figure}[h] 
   \centering
   \includegraphics[width=1.5in]{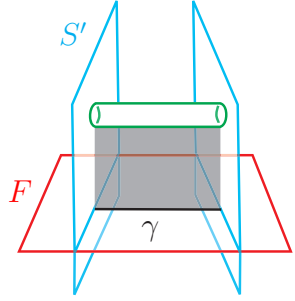} 
   \caption{Isotoping a stabilization (in green) into $X$ or $Y$ yields a boundary compressing disk based at $\gamma$.}
   \label{fig:stabilization_boundary_compressing_disk}
\end{figure}

\end{proof}

The bound in Corollary~\ref{cor:bound} can be weakened slightly to give a bound in terms of the number of components of $S \cap F$ instead of the number of components of $F$ cut along $S \cap F$ via the following lemma.

\begin{lemma}
\label{curves_of_intersection}
Let $F$ be a surface and let $\mathcal C \subset F$ be a disjoint union of $k$ simple closed curves. Then $F$ cut along $\mathcal C$ has at most $k + 1$ components. 
\end{lemma}

\begin{proof}
Let $c_1, \ldots, c_k$ be the components of $\mathcal C$. Now, $F$ cut along $c_1$ has either one or two components, depending on whether $c_1$ is non-separating or separating. Proceeding inductively, taking $F$ cut along $c_1 \cup \ldots \cup c_{j-1}$ and cutting along $c_j$ adds a new component if and only if $c_j$ is separating in $F$ cut along $c_1 \cup \ldots \cup c_{j-1}$. Thus, the extreme situation is if every component of $\mathcal C$ is separating in $F$ ({\em e.g.}~$\mathcal C$ consists entirely of inessential curves, like the resulting curves of intersection with $F$ when a Heegaard surface is isotoped to be close to a spine), in which case $F$ cut along $\mathcal C$ has $k + 1$ components. 
\end{proof}

Theorem~\ref{thm:main} now follows readily.

\begin{proof}[Proof of Theorem~\ref{thm:main}] 
Lemma~\ref{curves_of_intersection} implies that if \hs is a Heegaard splitting intersecting a mutually separating essential surface $F$ transversely in $k$ simple closed curves, then $|F \cap V| + |F \cap W| - 1 \leq k.$ The result now follows from Corollary~\ref{cor:bound}.
\end{proof}

The main theorem in \cite{DerbyTalbot2007} states that if $F$ is a torus and $V \cup_S W$ intersects $F$ in $2k$ simple closed curves essential on both surfaces, then $V \cup_S W$ is isotopic to an amalgamation along $F$ after at most $k$ stabilizations. This theorem can be seen to be a special case of Theorem~\ref{thm:main}, by first applying the techniques of Lemma 4.6 in \cite{DerbyTalbot2007} to construct an isotopy of $S$ with a surface that  intersects $F$ in $k$ inessential simple closed curves. 

\subsection{Non-separating surfaces and amalgamation}
\label{sec:nonseparating}

We conclude this section with a brief discussion of the situation where $F$ is not assumed to be mutually separating. We gave the example in Section~\ref{sec:amalgamation} that a connected, non-separating surface in $M$ is not mutually separating and thus cannot be used to form an amalgamation using Definition~\ref{def:amalgamation}. Every surface, however, can be made into a mutually separating one. 

Given a surface $F$ in $M$, form the surface $F'$ by adding to $F$ parallel copies of various components of $F$ until $F'$ becomes mutually separating. Note that for each component $\widetilde{F}$ of $F$ for which we add a parallel copy to form $F'$, the manifold $M$ cut along $F'$ has a component homeomorphic to $\widetilde{F} \times I$. Following the language of \cite{DerbyTalbot2006}, the surface $F'$ is called an {\em amalgamatable modification of $F$}. Theorem~\ref{thm:main} can thus be applied to a Heegaard splitting \hs and the surface $F'$.

Depending on how we form the surface $F'$ and how the Heegaard surface $S$ intersects $F'$, we may actually be able to use {\em fewer} stabilizations than the numbers prescribed by Corollary~\ref{cor:bound} and Theorem~\ref{thm:main} to obtain the requisite boundary compressing disks used in the proof of the corollary. As a general statement would be technical and cumbersome, we illustrate this idea by considering the special case mentioned before, that $F$ is connected and non-separating.

\begin{theorem}
\label{thm:surfacebundles}
Let $M$ be a 3-manifold containing a connected non-separating surface $F$. Assume that \hs is a Heegaard splitting such that $S$ intersects $F$ in $k$ inessential simple closed curves of intersection that bound disks in the same compression body (achieved, for example, by isotoping $S$ to be near a spine). Then \hs is isotopic to an amalgamation along two parallel copies of $F$ after at most $k+ 1 - \chi(F)$ stabilizations.
\end{theorem}

Note that $F$ by itself is not mutually separating in $M$, hence we need to add a parallel copy of $F$ to make a surface $F'$ to form an amalgamation. Theorem~\ref{thm:main} applied to $S$ and $F'$ would give nearly twice the bound obtained in the above theorem. 

\begin{figure}[h] 
   \centering
\includegraphics[width=5in]{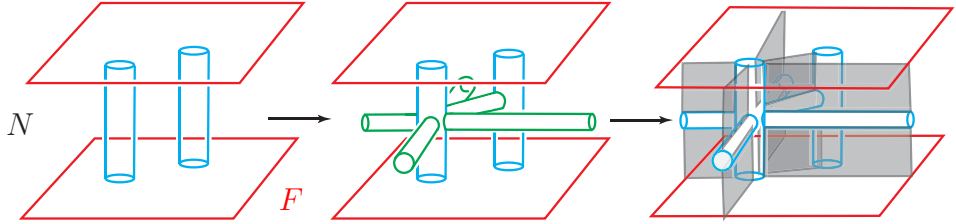} 
   \caption{Stabilizations (in green) of \hs in $F \times I$ give a complete set of boundary compressing disks for \underline{both} $F - \{ open \ disks\}$ and its parallel copy.}
   \label{fig:nonseparating}
\end{figure}

\begin{proof}
The proof is similar to the proof of Corollary~\ref{cor:bound}. Assume that $F'$ is obtained from $F$ and a parallel copy so that in the component $N$ homeomorphic to $F \times I$ between $F$ and its copy, $S \cap N$ is a disjoint union of vertical annuli. By assumption, $S$ intersects each component of $F'$ in inessential curves bounding disks in the same compression body. Without loss of generality assume that the curves bound disks in $V$. Let $\Gamma_W$ be a complete system of compressing arcs $\Gamma_W$ for $F \cap W$, and as in the proof of Corollary~\ref{cor:bound} stabilize \hs by adding 1-handles whose cores are the arcs in $\Gamma_W$. Pushing the 1-handles into $N$ yields a complete system of boundary compressing disks for {\em both} components of $F' \cap W$ (see Figure~\ref{fig:nonseparating}). Thus Theorem~\ref{theamalgamationlemma} implies the stabilized splitting is an amalgamation along $F'$, and the number of stabilizations added is $k +1 - \chi(F)$. 
\end{proof}

\psfrag{Ft}{\color{red}{$F$}}
\psfrag{N}{$N$}

\section{Bounding curves of intersection}
\label{sec:bounds}

\subsection{Upper bounds}

\begin{definition}
\label{def:minimal_intersection_number}
Let $S$ and $F$ be surfaces in a 3-manifold $M$. Define the {\em minimal intersection number} of $S$ and $F$ to be the minimum of $|S \cap F|$ over all possible isotopies of $S$ and $F$ where $S$ and $F$ intersect transversely.
\end{definition}

In light of Theorem~\ref{thm:main} it is of interest to bound above the minimal intersection number of a Heegaard surface $S$ and an essential surface $F$ to address the Stabilization Problem. There are examples where this can be done, for example Haken showed in \cite{Haken} that if $M$ is reducible, the minimal intersection number of $S$ with some set of essential spheres $F$ is bounded above by the number of components of $F$. Casson and Gordon generalized this in \cite{CassonGordon} to include the case that $F$ contains disks as well. 

If $F$ is connected and the genus of $F$ is one, then \cite{DerbyTalbot2007} implies that the minimal intersection number is bounded above by $4 g(S) - 4$, assuming that $F$ is a JSJ torus and $S$ can be isotoped to intersect $F$ in essential simple closed curves. If $F$ has genus greater than one, then Johannson has established upper bounds of $6g(S) - 11$ after possible modification of $F$ by annulus-compressions \cite{Johannson}, and more generally of $-48g(S)\chi(F)$ if $F$ is perhaps modified by Dehn twists along essential tori \cite{Johannson2}. These last three bounds are not necessarily strict.  Providing better upper bounds on the minimal intersection number is an avenue for further research. 

\subsection{Lower bounds}
We now consider the contrasting problem of finding lower bounds for the minimal intersection number of a Heegaard surface $S$ and an essential surface $F$. The surfaces $S$ and $F$ must intersect in at least $|F|$ simple closed curves. While there are several examples where this is also an upper bound (such as those given by Haken), Johannson writes in the introduction to his paper \cite{Johannson} that this upper bound ``cannot be expected to hold in general." Here, using a result of Schultens and Weidmann \cite{SchultensWeidmann}, we prove that Johannson's speculation is indeed true.  

\begin{definition}
Define the {\em Heegaard genus} of $M$, $h(M)$, to be the minimum of $g(S)$ where \hs is a Heegaard splitting of $M$. If $M$ is disconnected, then $h(M)$ is the sum of the Heegaard genera of the components of $M$.
\end{definition}

\begin{definition}
Let $F$ be a mutually separating essential surface in $M$. Define the {\em amalgamation genus of $M$ with respect to $F$}, $a(M,F)$, to be the minimal genus of a Heegaard splitting that is an amalgamation along $F$. 
\end{definition}

Note that if $M$ has boundary, then a Heegaard splitting of $M$ is really a Heegaard splitting of $M$ given some partition of the boundary components of $M$. In this case, we can further restrict the definitions of $h(M)$ and $a(M,F)$ by requiring that the minimal genus Heegaard splittings respect a given partition of the boundary components. 

\begin{remark}
\label{rmk:amalgamationgenus}
Suppose that $F$ is a mutually separating essential surface separating $M$ into $X$ and $Y$. Assume that $M$ is closed and that $F$ is connected. Then by amalgamating minimal genus Heegaard splittings of the components of $X$ and $Y$ we obtain the equation
\[ a(M,F) = h(X) + h(Y) - g(F). \]
\end{remark}

If $M$ is not closed or $F$ is disconnected, then this same equation only holds if, given a partition of $\partial M$, we take into account appropriate boundary partitions of $X$ and $Y$ in defining $h(X)$ and $h(Y)$. For example, the amalgamation of two ``Type I'' Heegaard splittings \cite{ScharlemannThompson} of two copies of $(surface) \times I$ along a single boundary component gives a minimal genus splitting of $(surface) \times I$, however the amalgamation of two ``Type II'' splittings does not (here, the boundary components of $(surface) \times I$ are partitioned into separate compression bodies in defining $a(M,F)$). 

Note that $h(M) \leq a(M,F)$ by definition. The situation where this is a strict inequality is of special interest.

\begin{definition}
\label{def:degeneration_of_Heegaard_genus}
Let $M$ be a 3-manifold containing an essential surface $F$ such that $h(M) < a(M,F)$. Then $M$ is said to have {\em degeneration of Heegaard genus}.
\end{definition}

See \cite{BachmanDT} for a general discussion of the notion of degeneration of Heegaard genus. Before proving Theorem~\ref{thm:lowerbound}, we establish a result about the minimal intersection number of Heegaard surfaces and essential surfaces in 3-manifolds that applies particularly in the case of degeneration of Heegaard genus.

\begin{theorem}
\label{thm:degeneration_of_Heegaard_genus}
Suppose that a 3-manifold $M$ contains an essential surface $F$, and let $d = a(M,F) - h(M)$. Then the minimal intersection number of a minimal genus Heegaard surface of $M$ and the essential surface $F$ is at least $d + \chi(F)$. 
\end{theorem}

\begin{proof}
Let \hs be a minimal genus Heegaard splitting of $M$ of genus $g = h(M)$, and assume that the minimal intersection number of $S$ and $F$ is $k$. By Theorem~\ref{thm:main}, \hs stabilizes to be an amalgamation along $F$ after at most $k - \chi(F)$ stabilizations. Hence, $$a(M,F) \leq g + k - \chi(F).$$ By assumption, $d = a(M,F) - h(M) = a(M,F) - g$, thus $a(M,F) = g + d$. Plugging this into the first inequality, we obtain $$g + d \leq g + k - \chi(F),$$ which implies $d + \chi(F) \leq k$. 
\end{proof}

We now prove Theorem~\ref{thm:lowerbound}.  
 
\begin{proof}[Proof of Theorem~\ref{thm:lowerbound}]
Let $n$ be a positive integer. By Theorem 27 in \cite{SchultensWeidmann}, there exist 3-manifolds $X_n$ and $Y_n$ with torus boundary (homeomorphic to a Seifert fibered space over a disk with $n$ exceptional fibers and a tunnel number $n$ knot complement, respectively) admitting minimal genus Heegaard splittings of genus $n$ and $n+1$, respectively, such that an amalgamation of these splittings is an $n$-times stabilization of a genus $n$ Heegaard splitting $V_n \cup_{S_n} W_n$ of $M_n = X_n \cup Y_n$. Thus $h(M_n) \leq n$, and by Remark~\ref{rmk:amalgamationgenus}, $a(M_n, T_n) = h(X_n) + h(Y_n) - 1 = 2n$, where $T_n = \partial X_n = \partial Y_n$ in $M_n$. Hence $d = a(M_n, T_n) - h(M_n) \geq 2n - n = n$. 

Assume that $S_n$ intersects $T_n$ in a minimal number $k$ of simple closed curves. By Theorem~\ref{thm:degeneration_of_Heegaard_genus}, $d \leq k$, which implies $n \leq k$. \end{proof}

\bibliographystyle{alpha}
\bibliography{070501_Stabilizations_amalgamations_intersectio_curves}

\end{document}